\documentstyle[12pt]{article} 
\newcommand{\doublespace}{
   \renewcommand{\baselinestretch}{1.2}
   \large\normalsize}

\def \Der{{\rm Der}}
\def \Z{\Bbb Z}
\def \C{\Bbb C}

\def \wt{{\rm wt}}

\def \End{{\rm End}}
\def \Hom{{\rm Hom}}

\def \<{\langle} 
\def \>{\rangle} 

\def \a{\alpha }

\def \l{\lambda }
\def \L{\Lambda }

\def \qed{\mbox{ $\square$}}
\def \pf {\noindent {\bf Proof:} \,}
\input amssym.def
\input amssym
\doublespace
\begin{document}
\newtheorem{thm}{Theorem}[section]
\newtheorem{prop}[thm]{Proposition}
\newtheorem{cor}[thm]{Corollary}
\newtheorem{lem}[thm]{Lemma}
\newtheorem{rem}[thm]{Remark}
\newtheorem{de}[thm]{Definition}
\begin{center}
{\Large {\bf Compact automorphism groups of vertex operator algebras}} \\
\vspace{0.5cm}
Chongying Dong\footnote{Supported by NSF grant 
DMS-9303374 and a research grant from the Committee on Research, UC Santa Cruz.}, Haisheng Li\footnote{Current address: Department of Mathematical Sciences,
Rutgers University, Camden, NJ, 08102}, Geoffrey Mason\footnote{Supported by NSF grant DMS-9122030 and a research grant from the Committee on Research, UC Santa Cruz.}
\\
Department of Mathematics, University of
California, Santa Cruz, CA 95064
\end{center}

\hspace{1.5 cm}
\begin{abstract} Let $V$ be a simple vertex operator algebra which admits
the continuous, faithful action of a compact Lie group $G$ of automorphisms.
We establish a Schur-Weyl type duality between the unitary, irreducible
modules for $G$ and the irreducible modules for $V^G$ which are contained
in $V$ where $V^G$ is the space of $G$-invariants of $V.$ 
We also prove a concomitant
Galois correspondence between vertex operator subalgebras of $V$ which contain
$V^G$ and closed Lie subgroups of $G$ in the case that $G$ is abelian. These
results extend those of [DM1] and [DM2].
\end{abstract}

\section{Introduction}

In this paper we study simple vertex operator algebras $V$ which admit
a continuous action of a compact Lie group $G$ (possibly finite). This
is a natural continuation of our earlier work [DM1], [DM2] in which we
limited ourselves to the case of finite solvable groups. Such pairs
$(V,G)$ occur naturally: for example the vertex operator algebras
$V_L$ associated with a positive definite even lattice $L$ (cf. [B],
[FLM]) admit Lie groups which arise from exponentiating the Lie
algebra associated with $L.$ A special case of this is when $L$ is
itself a root lattice (of type $ADE$) associated to a finite-dimensional
simple Lie algebra $\frak g$, in which case $V_L$ is the 
basic representation of the affine Lie algebra $\hat \frak g$ and $V_L$
admits a unitary action of the corresponding Lie group. This latter
situation was studied recently in a preprint of Baker and Tamanoi [BT],
in which they carried out some explicit calculations and made a general
conjecture (which we explain in due course). Our main theorem yields, as 
a very special case, a proof of their conjecture. It also yields an extension
of Theorem 3 of [DM1], originally proved only for solvable groups, to 
{\em all} finite groups.

The main result, Theorem \ref{t1} below, is a duality theorem of
Schur-Weyl type. Precisely, it says the following: let $I$ index the
finite dimensional irreducible (continuous) representations of $G,$
and for $\l\in I$ let $W_{\l}$ be the corresponding $G$-module. Then there is a decomposition of $V,$ considered as $G
\times V^G$-module, as follows:
$$(*)\ \ \ \ \ \ \ \ \ V=\bigoplus_{\l\in I}W_{\l}\otimes V_{\l}\ \ \ \ \ \ \ \  \ $$
where each $V_{\l}$ is a {\em non-zero, irreducible}
$V^G$-module. Here, $V^G$ is the vertex operator subalgebra of $V$
consisting of the $G$-invariants. Moreover, $V_{\l}$ and $V_{\mu}$ are
{\em inequivalent} $V^G$-modules if $\l\ne \mu.$ 
 Such results have been predicted --
at least for finite groups -- in the physics literature [DVVV]. The
decomposition $(*)$ was obtained for finite solvable groups in
[DM1]. The conjecture of [BT] is essentially the decomposition
$(*)$ in the case $V=V_L$ with $L$ of type $A_l$ and $G\simeq
SU(l+1).$ Note in particular that the Theorem implies that $V^G$ is
itself a simple vertex operator algebra, since $V^G$ is just the
subspace $V_1$ corresponding to the trivial $G$-module.
We also give another kind of duality or reciprocity result (Corollary
\ref{c2.7}) which compares the decomposition $(*)$ for 
a group $G$ and subgroup $H.$ Thus the pair $(G,V^G)$ behaves 
just like a classical dual pair [H2].

We also give an infinitesimal version of $(*)$ (Theorem
\ref{t2}).  Namely, we define the notion of a {\em derivation} of a
vertex operator algebra and observe that derivations form a Lie
algebra. Then the analogue of $(*)$ holds where $I$ is now the
set of dominant weights for a finite-dimensional semi-simple Lie
algebra of derivations of $V$ and $W_{\l}$ is the  corresponding
irreducible highest weight module.

As in [DM1], [DM2], $(*)$ also has applications to what we have
called {\em quantum Galois theory} (loc. cit.). Precisely, one is
looking for a Galois correspondence between  vertex operator
subalgebras of $V$ which contain $V^G$ and subgroups of $G.$ The case in
which $G$ is finite and either nilpotent or dihedral was discussed in
(loc. cit.), where the exact analogue of the usual Galois
correspondence was established. Here we consider the case of compact
abelian Lie groups (we could easily extend this to the case of compact
nilpotent Lie groups, along the line of Theorem 2 of [DM2]); for
infinite groups one expects only ``closed'' subgroups of $G$ to
intervene in the Galois correspondence, and we shall see that this is
so.  Namely, we prove that if $G$ is a compact abelian Lie group
acting continuously on the simple vertex operator algebra $V,$ then the map
$H\mapsto V^H$ gives a bijection between the {\em closed Lie subgroups}
of $G$ and the vertex operator subalgebras of $V$ which contain $V^G.$
For example, if $G=S^1$ is the circle then the vertex operator
subalgebras of $V$ which contain $V^{S^1}$ are, apart from $V^{S^1}$ itself,
precisely these of the form $V^F$ for a finite (cyclic) subgroup $F$ of
$S^1.$ Thus in this case the ``Krull topology'' one puts on $G$ is the fc
(closed=finite) topology.

The method of proof that we use is rather different from that of [DM1],
[DM2], though some key points are the same. We employ a certain $\Z$-graded
Lie algebra $\hat V$ naturally associated to $V$ that we have used elsewhere
[DLM], combined with some of Howe's ideas [H1] which we have adapted to
the present situation.

The main theorem and its infinitesimal analogue are proved in Section 2. The
Galois correspondence is proved in Section 3.

\section{The main theorem}
\setcounter{equation}{0}

Let $V=(V,Y,{\bf 1},\omega)$  be a simple vertex operator algebra 
and $G$ a group of
automorphisms of $V$ (cf. [FLM]). Assume that $G$ is a finite-dimensional 
compact Lie group and that the action of $G$ is
continuous. Let $\{W_{\l}|\l\in I\}$ be the set of inequivalent 
finite-dimensional irreducible representations. Then
any finite dimensional representation of $G$ is a direct sum of 
$W_{\l}$'s. Since each homogeneous subspace $V_n$ of $V$ is finite
dimensional, $V$ is a direct sum of finite-dimensional irreducible
$G$-modules:
$$\displaystyle{V=\bigoplus_{\lambda\in I}V^{\lambda}}$$
where $V^{\l}$ is the sum of all $G$-submodules of $V$ which are isomorphic
to $W_{\l}.$ In particular, the subspace of $G$-invariants $V^G$ is exactly
$V^{1}$ if $W_{1}$ is 
the 1-dimensional trivial representation. 

Consider the quotient space
$$\hat V={\C}[t,t^{-1}]\otimes V/D{\C}[t,t^{-1}]\otimes V$$
where $D={d\over dt}\otimes 1+1\otimes L(-1).$ Denote by $v(n)$ the image
of $v\otimes t^n$ in $\hat V$ for $v\in V$ and $n\in \Z.$ Then $\hat V$
is $\Z$-graded by defining the degree of $v(n)$ to be $\wt v-n-1$ if $v$
is homogeneous. Denote the homogeneous subspace of degree $n$ by $\hat V(n).$
The space $\hat V$ is, in fact, a $\Z$-graded Lie algebra with bracket
$$[a(m), b(n)]=\sum_{i=0}^{\infty}{m\choose i}a_ib(m+n-i)$$
(see [L1] and [DLM]). We define a $G$-action on $\hat V$ by $gv(n)=(gv)(n).$
Then $\hat V(n)$ is isomorphic to $V$ as $G$-modules. In particular,
$\hat V$ is a semi-simple $G$-module and $(\hat V)^G=\hat V^G.$
 The following proposition is proved in
[L1] and [DM1]:
\begin{prop}\label{p1} $V$ is spanned by $\{u_nv|u\in V,n\in\Z\}$ where 
$v\in V$ is a fixed nonzero vector. 
\end{prop}

For any $n\in \Z$ we define $V_{\leq n}=\sum_{m\leq n}V_m$ and 
$V_{>n}=\sum_{m>n}V_m.$ Then $V=V_{\leq n}\oplus V_{>n}$ and 
$\End V=\End V_{\leq n}\oplus \End V_{>n}\oplus \Hom_{\C}(V_{\leq n}, V_{>n})
\oplus \Hom_{\C}(V_{>n}, V_{\leq n}).$ Note that $\End V$ is a $G$-module
via conjugation and that $\End V_{\leq n}$ is a semi-simple submodule
as $V_{\leq n}$ is finite-dimensional. Let $T'$ denote
the projection of $T\in \End V$ onto $\End V_{\leq n}.$ Then this 
projection is a $G$-module homomorphism.   

Since $gY(u,z)g^{-1}=Y(gu,z)$ for $v\in V$ and $g\in G$ we see that
the map $v(m)\mapsto v_m$ is a $G$-module homomorphism from $\hat V$
to $\End V.$ Moreover, $f$: $v(n)\mapsto v_n'$ is a $G$-module
homomorphism from $\hat V$ to $\End V_{\leq n}$ for any $n$ as $g$
preserves both $V_{\leq n}$ and $V_{>n}.$ The 
following lemma is an immediate consequence of Proposition
\ref{p1}.
\begin{lem}\label{l1}
 The map $f$ is an epimorphism from the semi-simple $G$-module
$\hat V$ to the semi-simple $G$-module $\End V_{\leq n }$ for any $n.$
In particular, $f( \hat V^G)=(\End V_{\leq n})^G.$ 
\end{lem}

Fix $\l\in I.$ As the action of $G$ on $W_{\l}$ induces an epimorphism
$\C[G]\to \End W_{\l},$ we can define the space of highest weight vectors
in $W_{\l}$ as follows: fix a basis of $W_{\l},$ with $e_{ij}$ the 
corresponding standard basis of $\End W_{\l}.$ Then 
$w\in W_{\l}$ is a highest weight vector if $e_{ij}w=0$ whenever $j>i$ or
$i=j, i>1,$  and if $e_{11}w=w.$ We can then define the space $V_{\l}$ 
of highest weight vectors in $V^{\l}$ in the obvious way.

\begin{lem}\label{l2} Each $V_{\l}$ is a module for $V^G.$
\end{lem}

\pf Since the action of $G$ commutes with the action of vertex
operator $Y(u,z)$ for $u\in V^G,$ the lemma follows immediately.
\qed

Note that as both  $G$-module and $V^G$-module,
$V^{\l}$ and $W_{\l}\otimes V_{\l}$ are isomorphic. Here,
$G$ acts on the first tensor factor of $W_{\l}\otimes V_{\l}$ 
and $V^G$ acts on the second tensor factor.
The main result of this paper is the following.
\begin{thm}\label{t1} (1) Each $V_{\l}$ is an irreducible $V^G$-module.

(2) $V^{\l}$ is non-zero for any $\l\in I.$   

(3) $V_{\l}$ and $V_{\mu}$ are isomorphic $V^G$-modules if, and only if,
$\l=\mu.$

\end{thm}

\pf The main idea in the proof of (1) and (3) goes back to Howe's theory
of dual pairs (see Theorem 8 of [H1]).

For (1) it is sufficient to show that as a $V^G$-module, $V_{\l}$ is
generated by any nonzero vector of $V_{\l}.$  
First choose $x,y\in V_{\l}$ to be homogeneous and 
linearly independent and $n\in \Z$ such that $x,y\in V_{\leq n}.$
Then as a $G$-module we can write $V_{\leq n}$ as a direct sum of
three submodules $W_{\l}\otimes x\oplus W_{\l}\otimes y\oplus W.$
Define a map $\alpha\in \End V_{\leq n}$ such that $\a(u\otimes
x+v\otimes y+w)=v\otimes x+u\otimes y+w$ for $u,v\in W_{\l}$ and $w\in
W.$ Then it is easy to see that $g\alpha g^{-1}=\alpha$ for all $g\in G.$
That is, $\alpha\in (\End V)^G.$ 

On the other hand by Lemma \ref{l1} there exist homogeneous 
vectors $v^1,...,v^k\in V^G$ and $m_1,...,m_k\in\Z$ such that 
$(\sum_{i=1}^kv^i_{m_i})'=\a.$ Thus $(\sum_{i=1}^kv^i_{m_i})'x=y.$
Since $x,y$ and $v_i$ are homogeneous there exists a subset $J$
of $\{1,...,k\}$ such that $\sum_{i\in J}v^i_{m_i}x=y.$ That is, $V_{\l}$ is
generated by any nonzero homogeneous vector of $V_{\l}.$ 
If $x$ is not homogeneous we can use the operator $L(0)$ to produce
a nonzero homogeneous vector $x'\in V_{\l}$ and the argument above 
gives (1).

For (2) we first assert  that there exists a faithful $G$-module $W$
such that $W=W_{\l_1}+\cdots+W_{\l_k}$ with $V^{\l_i}\ne 0.$ 
Set $K_n=\ker (G, V_{\leq n}).$ Then $\cdots K_n\supset K_{n+1}\cdots$ is 
 a descending chain of closed 
subgroups of $G.$ Since $G$ is a finite-dimensional Lie group, there exists  $N\geq 0$ such that $\dim K_n$
is a constant if $n\geq N.$ Let $K_n^0$ be the connected component
of $K_n$ which contains the identity $e.$ Then $K_n^0/K_n^0\cap K_{n+1}
\simeq K_n^0K_{n+1}/K_{n+1}\subset K_n/K_{n+1}.$ As $K_n^0/K_n^0\cap K_{n+1}$
is connected and compact and $K_n/K_{n+1}$ is a finite group if
$n\geq N,$  $K_n^0/K_n^0\cap K_{n+1}$ is a one point set. That is, 
$K_n^0=K_{n+1}^0.$ This implies that $K_N^0$ acts trivially on $V.$ Thus
$K_n^0=e$ if $n\geq N$ and $K_n$ is  finite. Consequently,
$K_n=e$ and the representation of $G$ on 
$V_{\leq n}$ is 
faithful if $n$ is big enough. For such $n$ we may take  $W=W_{\l_1}+\cdots+W_{\l_k}$ where $W_{\l_1},...,W_{\l_k}$ are all the distinct $W_{\l}$'s
contained in $V_{\leq n}.$ Our assertion follows.

Let $\Lambda=\{\l\in I|V^{\l}\ne 0\}.$ It was proved in [DM1] that 
if $G$ is a finite group and $\l_1,\l_2\in \L$ and $W_{\gamma} $ for 
$\gamma\in I$  is a submodule
of $W_{\l_1}\otimes W_{\l_2}$ then $\gamma\in \L$ (see the proof
of Theorem 2 of [DM1]). The same result holds here with the same proof.
We refer the reader to [DM1] for details. 

Next we show that if $\l\in \L$ then $\l^*\in \L$ where $W_{\l^*}$ is
the conjugate representation. By Proposition \ref{p1} $V$ is spanned
by $\{u_nW_{\l}|u\in V, n\in\Z\}.$ For fix $n\in \Z$ and $\mu\in \L$
observe that the span of $\{u_nW_{\l}|u\in W_{\mu}\}$ is a
$G$-submodule of $W_{\l}\otimes W_{\mu}.$ Note that the trivial
$G$-module necessarily occurs in $V,$ and that the trivial module is a
submodule of $W_{\l}\otimes W_{\mu}$ if, and only if, $\nu=\l^*.$ So
indeed we have $\l^*\in \Lambda$ whenever $\l\in\Lambda.$

It is well-known that every irreducible module is contained in
$M^{\otimes m}\otimes \bar M^{\otimes n}$ for some positive integers
$m,n$ whenever $M$ is a faithful 
representation of $G$ and $\bar M$ is its conjugate (eg, Theorem 4.4
of [BD]). Now $W$ is a faithful representation of $G$ inside $V$ and its
conjugate
$\bar W =\sum W_{\l^*_i}$ also lies in $V.$ Thus $\L=I.$ This proves (2).

The proof of (3) is similar to that of (1).  We first take $n$
such that $V^{\l}_n\ne 0,$ as we may by (2).
 Then $V_{\le n}$ is a direct sum of three
$G$-modules $V^{\l}_n\oplus V^{\mu}_n\oplus M$ for a suitable
$G$-submodule $M$ of $V_{\le n}.$ Define $\beta\in \End V_{\le n}$
such that it is the identity on $V^{\l}_n$ and zero on $V^{\mu}_n$ and $M.$
Clearly $\beta\in ( \End V_{\le n})^G.$ Use Lemma \ref{l1} to obtain
homogeneous vectors $v^1,...,v^k$ of $V^G$ such that $\sum v^i_{\wt v^i-1}$
is the identity on $V^{\l}_n$ and zero on  $V^{\mu}_n$ and $M.$ Thus there is
no $V^G$-module homomorphism between $V_{\l}$ and $V_{\mu},$  as required.
\qed. 

\begin{cor}\label{c1} Let $V$ be a simple vertex operator algebra and $G$ be a
finite automorphism group of $V.$ Then 
$$V=\sum_{\chi\in\hat G}W_{\chi}\otimes V_{\chi}$$
where $W_{\chi}$ is an irreducible $G$-module affording $\chi$ and
$V_{\chi}$ is an irreducible $V^G$-module. Moreover, $V_{\chi}$
and $V_{\psi}$ are non-zero inequivalent $V^G$-modules if $\chi$ and $\psi$ are
different irreducible characters of $G.$
\end{cor}

\begin{rem} Corollary  \ref{c1} was proved in [DM1] in the case that $G$ is
a finite solvable group (loc. cit., Theorems 2 and 3). 
\end{rem}

There is a reciprocity law attached to the dual pair $(G,V^G)$ as in the 
classical
case (cf. [H2]). Let $H$ be a compact Lie subgroup of $G.$ Then by Theorem
\ref{t1} there is a decomposition of $V$ into irreducible 
$H\times V^H$-modules
$$ V=\bigoplus_{\mu\in J}X_{\mu}\otimes V_{\mu}$$
where $J$ indexes
 the finite dimensional irreducible (continuous) representations
of $H,$  $X_{\mu}$ is the corresponding $H$-module for
$\mu\in J,$ and $V_{\mu}$ is the corresponding irreducible $V^H$-module.

Let us decompose the representations of $G$ into a direct sum of $H$-modules:
$$W_{\l}=\bigoplus_{\mu\in J} m_{\l\mu}X_{\mu}$$
for $\l\in I$ where $m_{\l\mu}$ is the multiplicity of $\mu$ occurring in $\l.$
By Theorem \ref{t1}, $V$ is a complete reducible $V^G$-module.
Theorefore, the $V^H$-modules $V_{\mu}$ for $\mu\in J$ decompose into 
a direct sum of irreducible $V^G$-modules:
$$V_{\mu}=\bigoplus_{\l\in I}n_{\mu\l}V_{\l}.$$
Note that both $m_{\l\mu}=\dim \Hom_{H}(W_{\l},X_{\mu})$ and 
$n_{\mu\l}=\dim \Hom_{V^G}(V_{\mu},V_{\l})$ are finite.
\begin{cor}\label{c2.7} The following reciprocity holds:
$$m_{\l\mu}=n_{\mu\l}$$
for all $\l\in I$ and $\mu\in J.$
\end{cor}

\pf The proof is exactly the same as that in Section 5.7 of [H2] in the
classical case. For completeness we repeat this short proof. 

Note that $V$ is a complete reducible $H\times V^G$-module and consider
the decomposition of $V$ as a $H\times V^G$-module. Let $\l\in I$
and $\mu\in J.$ Clearly, in the decomposition
of $V$ into $H\times V^H$-modules, $X_{\mu}\otimes V_{\l}$ only occurs in
$X_{\mu}\otimes V_{\mu}$ with multiplicity $n_{\mu\l}.$  Similarly,
in the decomposition
of $V$ into $G\times V^G$-modules,  $X_{\mu}\otimes V_{\l}$ only occurs in
$W_{\l}\otimes V_{\l}$ with multiplicity $n_{\mu\l}.$ Thus
$m_{\l\mu}=n_{\mu\l},$ as claimed.   \qed

Next we present a Lie algebra version of Theorem \ref{t1}. We need the 
concept of a derivation of a vertex operator algebra. 
A derivation $D$ of vertex operator algebra $V$ is a degree-preserving 
linear map of $V$ such that $D{\bf 1}=D{\omega}=0$ and
$[D,Y(u,z)]=Y(Du,z)$ for all $u\in V.$ Denote the set of all derivations
of $V$ by $\Der V.$ Then $\Der V$ is a Lie algebra. Let
${\frak g}$ be a semi-simple Lie subalgebra of $\Der V.$ Then $V$ is 
completely reducible ${\frak g}$-module: $V=\oplus_{\lambda\in P}V^{\lambda}$
where $P$ is the set of dominant weights of ${\frak g}$ and
$V^{\l}$ is the sum of all irreducible ${\frak g}$-modules which are isomorphic
to the highest weight module $L(\l)$ 
for ${\frak g}$ with highest weight $\l.$
Let $V_{\l}$ be the space of all highest weight vectors in $V^{\l}.$
Then $V^0=V_0=V^{\frak g}$ is the space of ${\frak g}$-invariants 
of $V.$ That is, $V^{\frak g}=\{v\in V|{\frak g}v=0\}.$ It is easy to show 
that each $V_{\l}$ is a module for $V^{\frak g}.$ Then
$V^{\l}=L(\l)\otimes V_{\l}$ as ${\frak g}\otimes V$-modules where
${\frak g}$ acts on the first tensor factor and $V^{\frak g}$ acts on the
second tensor factor. Now use the proof of Theorem \ref{t1} to get
\begin{thm}\label{t2} $V^{\l}$ is non-zero for all $\l\in P,$ and 
$V_{\l}$ is an irreducible $V^{\frak g}$-module. Moreover, $V_{\l}$ 
and $V_{\mu}$ are isomorphic $V^{\frak g }$-modules if, and only, if
$\l=\mu.$ In particular,
$V^{\frak g}$ is a simple vertex operator algebra.
\end{thm}

Now we apply this result to the vertex operator algebras associated 
to highest weight representations of affine Lie algebras. Let ${\frak g}$ 
be a finite-dimensional simple Lie algebra and $\hat {\frak g}$ be
the corresponding affine Lie algebra.   Denote
by $L(\ell,\lambda)$ the highest weight $\hat{{\frak g}}$-module of level 
$\ell$ (which is a fixed positive integer) with highest weight $\lambda\in P.$ 
It is well known (cf. [DL], [FZ], [L2]) that
$L(\ell,0)$ is a simple vertex operator algebra. It is clear that ${\frak g}$
acts on $L(\ell,0)$ as derivations. The following corollary is an immediate
consequence of Theorem \ref{t2}.

\begin{cor} $L(\ell,0)$ has the following decomposition:
$$L(\ell,0)=\bigoplus_{\l\in P}L(\l)\otimes L(\ell,0)_{\l}$$
where each $L(\ell,0)_{\l}$ is an irreducible $L(\ell, 0)^{\frak g}$-module
and $L(\ell,0)_{\l}$ and $L(\ell,0)_{\mu}$ are isomorphic $L(\ell, 0)^{\frak g}$-modules if, and only, if $\l=\mu.$ 
\end{cor}

\section{Galois correspondence}

We will prove
\begin{thm}\label{t3.1} Let $V$ be a simple vertex operator algebra and let $G$ be a compact
abelian Lie group of (faithful , continuous) automorphisms of $V.$ Then the
correspondence
$$H\to V^H$$
sets up a bijection between the closed Lie subgroups of $G$ and the
vertex operator subalgebras of $V$ which contain $V^G.$ 
\end{thm}

With the main Theorem \ref{t1} now available, the proof follows the lines
of the corresponding result (Theorem 1 of [DM1]) for finite abelian groups.

We know [BD] that $G$ is isomorphic to $A\times T^n$ where $A$ is a
finite abelian group and $T^n$ the $n$-torus. And Pontryagin duality
says that the dual group $\hat G$ of unitary characters is isomorphic
to $\hat A\times \Z^n.$ The closed Lie subgroups of $G$ are precisely
the abstract subgroups of the form $\hat K$ where $K$ is a subgroup of 
$\hat G$ and $\hat K=\{g\in G|\l(g)=1 \  for \  all\ \l\in K\}.$

Now in the present case, Theorem \ref{t1} says the following: V has a
decomposition into distinct (non-zero) irreducible $V^G$-modules of the 
form
$$V=\bigoplus_{\l\in\hat G} V^{\l}.$$
Then for any subgroup $K$ of $\hat G$ we may define
$$V_K=\bigoplus_{\l\in K} V^{\l}.$$
It is clear that $V_K=V^{\hat K},$ and that $V_K=V_L$ if, and only if,
$K=L.$ So it suffices to show that each $V_K$ is a vertex operator 
subalgebra of $V$ and that every vertex operator subalgebra of $V$ containing
$V^G$ has this form.

Now if $\l_1,\l_2\in K< \hat G$ and if $u\in V^{\l_1},$ $v\in V^{\l_2},$ then
certainly $Y(u,z)v$ has component vectors in $V^{\l_1+\l_2}.$
So $V_K$ is certainly a vertex operator subalgebra of $V.$ On the 
other hand if $V^G\subset W\subset V$ and $W$ is a vertex operator subalgebra
then $W=\bigoplus_{\sigma\in S}V^{\sigma}$ for some subset $S\subset \hat G.$
Now as in the proof of Theorem \ref{t1} above (see also the proof of Theorem
2 of [DM1]) we see that $S$ is indeed closed under multiplication  and 
taking inverse and hence is a subgroup of $\hat G.$ This completes the
proof of Theorem \ref{t3.1}.

\end{document}